\documentclass[10pt]{article}

\usepackage[top=1in, bottom=1in, left=1.25in, right=1.25in]{geometry}

\usepackage[english]{babel}
\usepackage{amsmath}
\usepackage{graphicx}
\usepackage{amssymb}

\usepackage{algorithm}
\usepackage{algorithmic}

\usepackage{amsthm}

\newcommand{\norm}[2]{\lVert #1 \rVert_{#2}}
\newcommand{\normtwo}[1]{\lVert #1 \rVert_2}
\newcommand{\normB}[1]{\lVert #1 \rVert_B}
\newcommand{\define}{\stackrel{\text{def}}{=}}
\newcommand{\prior}{\Gamma_\text{prior}}
\newcommand{\noise}{\Gamma_\text{noise}}
\newcommand{\bigO}{{\cal{O}}}
\newcommand{\post}{\Gamma_\text{post}}

\newcommand{\red}{H_\text{red}}

\newcommand{\bx}{\textbf{x}}
\newcommand{\at}[2]{\left.#1\right|_{#2}}
\newcommand{\normal}{ {\cal{N}}}

\newcommand{\logdet}{\log\text{det}  }
\newcommand{\trace}{\text{Trace}}

\title{Fast computation of Uncertainty Quantification measures in the Geostatistical approach to solve Inverse Problems}
\author{Arvind K. Saibaba  \thanks{Department of Electrical and Computer Engineering, Tufts University,  196 Boston Avenue, Medford, MA - 02155 arvind.saibaba@tufts.edu}\and Peter K. Kitanidis  \footnotemark[1] \thanks{Department of Civil and Environmental Engineering, Yang \& Yamazaki Environment \& Energy Building  473 Via Ortega, Stanford, CA 94305 and Institute for Computational and Mathematical Engineering,  Huang Engineering Center
475 Via Ortega, Stanford, CA 94305, peterk@stanford.edu}}

\begin{document}
\maketitle

\begin{abstract}
We consider the computational challenges associated with uncertainty quantification involved in parameter estimation such as seismic slowness and hydraulic transmissivity fields. The reconstruction of these parameters can be mathematically described as Inverse Problems which we tackle using the Geostatistical approach. The quantification of uncertainty in the Geostatistical approach involves computing the posterior covariance matrix which is prohibitively expensive to fully compute and store. We consider an efficient representation of the posterior covariance matrix at the maximum a posteriori (MAP) point as the sum of the prior covariance matrix and a low-rank update that contains information from the dominant generalized eigenmodes of the data misfit part of the Hessian and the inverse covariance matrix. The rank of the low-rank update is typically independent of the dimension of the unknown parameter. The cost of our method scales as $\bigO(m\log m)$ where $m $ dimension of unknown parameter vector space. Furthermore, we show how to efficiently compute  measures of uncertainty that are based on scalar functions of the  posterior covariance matrix. 
The performance of our algorithms is demonstrated by application to model problems in synthetic travel-time tomography and steady-state hydraulic tomography.  We explore the accuracy of the posterior covariance on different experimental parameters and show that the cost of approximating the posterior covariance matrix depends  on the problem size and is not sensitive to other experimental parameters.
\end{abstract}

\section{Introduction}
One of the central challenges in the field of geosciences is to develop computationally efficient statistical methods for optimizing the use of limited and noisy environmental data to accurately estimate heterogeneous subsurface geological properties. In addition, it is necessary to quantify the corresponding predictive uncertainty. Mathematically, imaging can be performed using inverse problems theory, which uses measurements to make inference of system parameters. Efficient algorithms for inverse problems are necessary to solve problems of realistic sizes, quantified by  the spatial resolution of the reconstructed parameters and number of measurements available for reconstruction. Using these efficient algorithms scientists can gain better knowledge of soil moisture content, the porosity of geologic formations, distributions of dissolved pollutants, and the locations of oil deposits or buried liquid contaminants. These detailed images can then be used to better locate natural resources, treat pollution, and monitor underground networks associated with geothermal plants, nuclear waste repositories, and carbon dioxide sequestration sites. We aim to solve these problems by employing the geostatistical approach that stochastically models unknowns as random fields and uses Bayes' rule to infer unknown parameters by conditioning on measurements. However, due to high computational costs in identifying small scale features, these methods are challenging. These costs occur because solving inverse problems requires multiple expensive simulations of partial differential equations as well as representing high dimensional random fields, especially on irregular grids and complicated domains. Additional details about the Geostatistical approach have been provided in Section~\ref{sec:geostatistical}.

Uncertainty  in the context of Bayesian inverse problems is represented by the posterior probability density function. For linear inverse problems, if the measurement noise is additive Gaussian and the prior model is specified by a Gaussian random field, then the resulting posterior probability density function (PDF) is also Gaussian and is fully specified by calculating the mean (which coincides with the maximum a posteriori, or MAP, estimate), and the posterior covariance matrix. Computing the mean leads to a weighted linear least squares optimization problem, which can be tackled by several efficient numerical algorithms. For nonlinear inverse problems, a linearization of the measurement operator yields a local Gaussian for the posterior PDF. The MAP point can be computed by solving a weighted regularized nonlinear least squares problem and the posterior covariance matrix can be approximated by the inverse of the Hessian of the negative log-likelihood of the posterior PDF computed at the MAP estimate. Although considerable effort has been devoted to computing the MAP estimate (see for example~\cite{kitanidis1995quasi,saibaba2012efficient,saibaba2012application}), relatively fewer number of works have addressed the computation of the posterior covariance matrix. In this work we focus on an efficient representation of the posterior covariance matrix which is a measure of uncertainty associated with the reconstruction of the parameters of interest.

Computing and storing the approximation to the posterior covariance matrix is computationally infeasible because  the prior covariance matrices arising from finely discretized fields and certain covariance kernels are dense, and computing the dense measurement operator requires solving many forward PDE problems, which can be computationally intractable. In ill-posed inverse problems, the data is informative only about a low-dimensional manifold in the parameter space. This property has been exploited previously in developing efficient approximate representations to the posterior covariance matrix as the sum of the prior covariance matrix and a low-rank update that contains combined information from both the prior and the data misfit part of the Hessian (see for example,~\cite{bui2012extreme,flath2011fast,bui2013computational}). The low-rank modification is computed by the solution of a large-scale eigenvalue problem involving the prior-preconditioned Hessian. The prior covariance matrices can be modeled as discrete representations of operators of the form $\mathcal{A}^{-\alpha}$, where $\mathcal{A}$ is a partial differential operator (for e.g., the Laplacian) and $\alpha $ is a parameter chosen such that the infinite dimensional formulation is well-posed~\cite{stuart2010inverse}. Another choice for prior covariance matrices  is using Spartan Gibbs random field~\cite{hristopulos2003spartan}. In this work we focus on the Mat\'{e}rn class of covariance kernels~\cite{rasmussen2006gaussian}.

The ability of being able to compute measures of uncertainty is extremely important for the field of Optimal Experimental Design (OED), which seeks to determine the experimental setups which maximize the amount of information that can be gained about the parameters of interest. The design variables which control the accuracy of the parameter reconstructions could be the measurements or measurement types, numbers, locations of sources and/or detectors and other experimental conditions. A prevalent approach to OED involves optimizing an objective function which involves a scalar measure of uncertainty associated with the parameter reconstruction (defined on the basis of the posterior covariance matrix) and attempts to minimize this objective function with respect to design parameters. Since during the context of optimizing the experimental setup, the inverse problem has to be solved several times and the resulting uncertainty needs to be estimated at each iteration of the optimization routine, we would like an efficient method for computing the objective function (i.e., the measure of uncertainty).  We will provide an efficient method for computing a few of these uncertainty measures.  A good review of optimal experimental design in the Bayesian context is provided in~\cite{chaloner1995bayesian}. Common optimality criteria which can be used as objective functions include the alphabetic criteria, for example, A-, C-, D-, E-, and T- optimality criteria (these will be defined in Section~\ref{sec:uncertmeas}).  The definitions of the optimality criteria in the geostatistical context, along with a discussion of physical and statistical significance of these criteria and its applicability in non-Gaussian settings is available in~\cite{nowak2010measures}.

\textbf{Contributions}: 
We model the prior covariance matrix $\prior$ with entries arising from covariance kernels, as is common practice~\cite{kitanidis1995quasi}. Although the resulting covariance matrices are dense, in our previous work~\cite{saibaba2012efficient,ambikasaran2013large}, we have shown that we can obtain the best estimate using techniques (such as FFT based methods and $\mathcal{H}$-matrix approach) to reduce the storage and computational cost from $\bigO(m^2)$ to $\bigO(m\log m)$ (m is the number of unknown parameters). However, except when the number of measurements are small, e.g., $\bigO(100)$, computing entries of the posterior covariance matrix is computationally impractical. We show how to compute an efficient representation of the posterior covariance matrix as a low-rank modification of $\prior$ and the low-rank update is computed efficiently using a randomized algorithm. A major advantage of our approach is that there is great flexibility in experimenting with several covariance kernels, since the prior covariance matrix computations are handled in a black-box fashion. This is a major difference in our work compared to~\cite{bui2012extreme,flath2011fast,bui2013computational}, that we consider directly modeling the entries of the prior covariance matrix using the Mat\'{e}rn class of covariance kernels, instead of modeling the prior as the inverse of a discretized  differential operator (such as the Laplacian). Although these two approaches appear disparate, their equivalence has been established in~\cite{lindgren2011explicit}.    

A second contribution of this paper is that we provide an algorithm for approximating the posterior covariance matrix  that does not require forming the square root (or equivalently Cholesky factorization) and inverse of the prior covariance matrix $\prior$. The approach in ~\cite{bui2012extreme,flath2011fast,bui2013computational} considers the prior preconditioned Hessian (defined as $\prior^{1/2}\red\prior^{1/2}$, where $\red$ is the Gauss-Newton Hessian of the data misfit term. Computing the  square root of a matrix is an expensive operation for finely discretized grids arising from large-scale 3D problems. The work in~\cite{bui2012extreme,flath2011fast,bui2013computational} avoids this issue by considering priors for which the square-root is explicitly known. Since the matrix square root is not explicitly known for arbitrary covariance matrices, this assumption is very restrictive from a modeling stand point. The algorithm we propose only requires forming matrix-vector products (henceforth, referred to as matvecs). The key idea is to consider an equivalent generalized eigenvalue problem $Ax=\lambda Bx$ where $A $ is the Hessian corresponding to the data-misfit and $B = \prior^{-1}$ is the inverse prior covariance matrix. The randomized algorithm that we propose for approximating the posterior covariance matrix $\post$ is simple to implement, is computationally efficient, and comes with error bounds established in our previous work~\cite{saibaba2013randomized}.

Another important contribution of our work is the efficient computation of various measures of uncertainty which leverages the efficient representation of the posterior covariance matrix, written as a low-rank correction to the prior covariance matrix. 
A second computational burden occurs when the number of measurements are large because of operations on a dense cross-covariance matrix, which scale as $\bigO(n^3)$, where $n$ is the number of measurements. While some of the criteria (A- and C-) can be evaluated when the number of measurements are small, other criteria (such as D- and E-) are altogether computationally infeasible. However, using an efficient representation of the approximate posterior covariance and using matrix-free techniques, we show how several of these optimality criteria can be computed more efficiently. We note a further advantage of using covariance kernels to model $\prior$: Computing the variance of the posterior covariance requires computing the diagonals of the prior covariance matrix which can be easily computed, when the covariance kernel is specified explicitly. The application of computing these uncertainty measures in the context of optimal experimental design has been also covered in~\cite{alexanderian2014bayesian,alexanderian2014fast,alexanderian2014optimal}.

Finally, the efficiency of our proposed algorithms is demonstrated through challenging applications of estimating seismic slowness using traveltime ray tomography  (Section~\ref{sec:raytomog}) and estimating hydraulic conductivity (or transmissivity) using steady state hydraulic tomography  (Section~\ref{sec:ssht}).

\section{Geostatistical approach to inverse problems}\label{sec:geostatistical}

The Geostatistical approach (described in the following papers~\cite{kitanidis1995quasi,kitanidis2010bayesian,kitanidis2007on}) is one of the prevalent approaches to solve inverse problems. The geostatistical approach represents the unknown field to be estimated as a random field composed as the sum of a few deterministic term, typically low order polynomials, and a stochastic term which models small scale variability. The structure of the stochastic term is represented through the prior probability density function, which in practical applications is often parameterized through variograms and generalized covariance functions. The method has found several applications because it is generally practical and has the ability to quantify uncertainty. The method can generate best estimates, which can be determined in a Bayesian framework as a posteriori mean values or most probable values, measures of uncertainty as posterior variances or credibility intervals, and conditional realizations, which are sample functions from the ensemble of the posterior probability distribution.

In this section, we briefly review the geostatistical approach.  Let $s(\bx)$, the function to be estimated, be  modeled by a Gaussian random field. After discretization, we can write $s\sim \normal (X\beta, \prior) $ where $X \in \mathbb{R}^{m\times p}$ is a matrix of low-order polynomials known as the drift matrix, $\beta \in \mathbb{R}^p$ are a set of drift coefficients to be determined and $\prior$ is a covariance matrix with entries $\prior(i,j) = \kappa(\bx_i,\bx_j)$, and $\kappa(\cdot,\cdot)$ is a generalized covariance kernel~\cite{christakos1984problem,matheron1973intrinsic}. We have that $s \in \mathbb{R}^m$ and $m$ refers to the grid size.  In several instances, the measurements are corrupted by noise and we model these as Gaussian random variables. The measurement equation can be written as 
\begin{equation}
 y = h(s) + v \qquad v \sim \normal (0,\noise)
\end{equation}
where $y \in \mathbb{R}^{n}$ represents the noisy measurements, $h:\mathbb{R}^{m}\rightarrow\mathbb{R}^{n}$ is known as the {\it measurement operator} or {\it parameter-to-observation map} and $v$ is a random vector of observation error with mean zero and covariance matrix $\noise$. The matrices $\noise$, $\prior$ and $X$ are part of a modeling choice and more details to choose them can be obtained from the following reference~\cite{kitanidis1995quasi}. The parameters to be reconstructed are the values of the function at the grid locations $s$ and $\beta$ which are the drift coefficients.

We can write down the following expressions for the probability density functions (PDF) for the prior and the measurements 

\begin{equation}
\label{eqn:priorpdf}
p(s|\beta) \propto \exp\left( -\frac{1}{2}\norm{s - X\beta}{\prior^{-1}}^2\right)
\end{equation}

and 
\begin{equation}
\label{eqn:measurementpdf}
p(y|s) \propto \exp\left( -\frac{1}{2}\norm{y - h(s)}{\noise^{-1}}^2\right)
 \end{equation}
where we have used the vector norm $\|x\|_M = \sqrt{x^TMx}$ and $M$ is defined for an arbitrary symmetric positive definite matrix. Inference of the parameters from the measurements is obtained by invoking the Bayes' theorem, through the posterior PDF which is the product of two parts - the likelihood of the measurements and the prior distribution of the parameters. Employing Bayes' rule and assuming that the prior PDF for $\beta$ is uniform, i.e. $p(\beta) \propto 1$, we can write the expression for the PDF of the posterior distribution as 
\begin{equation}
\label{eqn:posteriorpdf}
p(s,\beta|y) \propto \exp\left( -\frac{1}{2}\norm{y - h(s)}{\noise^{-1}}^2 -\frac{1}{2}\norm{s - X\beta}{\prior^{-1}}^2\right) 
\end{equation}
The MAP estimate is computed by  maximizing the negative log likelihood of the posterior displayed in Equation~\eqref{eqn:posteriorpdf}, and can be obtained by solving the following optimization problem

\begin{equation}
\label{eqn:mapestimate}
\underset{\hat{s},\hat{\beta}}{\operatorname{arg\text{ }min}}\quad \frac{1}{2}\norm{y - h(s)}{\noise^{-1}}^2 + \frac{1}{2}\norm{s - X\beta}{\prior^{-1}}^2
\end{equation}

The choice of prior covariance kernels is the Mat\'{e}rn family of covariance kernels~\cite{rasmussen2006gaussian} 
\begin{equation}\label{eqn:maternfamily}
\kappa(\bx,\textbf{y}) = C_{\alpha,\nu}(r) = \frac{1}{2^{\nu-1}\Gamma(\nu)} (\sqrt{2\nu}\alpha r)^\nu K_\nu(\sqrt{2\nu}\alpha r) 
\end{equation}
where $r=\normtwo{\bx-\textbf{y}}$, $\Gamma$ is the Gamma function, $\alpha$ is a scaling factor, $K_\nu(\cdot)$ is the modified Bessel function of the second kind of order $\nu$. Equation~\eqref{eqn:maternfamily} takes special forms for certain parameters $\nu$. For example, when $\nu=1/2$, $C_{\alpha,\nu}$ corresponds to the exponential covariance function when $\nu = 1/2+\tilde{n}$ and	 $\tilde{n}$ is an integer, $C_{\alpha,\nu}$ is the product of an exponential covariance and a polynomial of order $\tilde{n}$. In the limit as $\nu\rightarrow\infty$, and for appropriate scaling of $\alpha$, $C_{\alpha,\nu}$ converges to the Gaussian covariance kernel. For a more detailed discussion of permissible covariance kernels, we refer the reader to the following references~\cite{matheron1973intrinsic,christakos1984problem}. Other possible choices for modeling the unknown fields are using the Spartan Gibbs random field models~\cite{hristopulos2003spartan}. Estimating the covariance parameters can be accomplished using the restricted maximum likelihood approach that has been outlined in~\cite{kitanidis1995quasi}. The choice of the parameter $\nu$ depends on the a priori information available of the smoothness of the field we wish to  reconstruct. The effect of the parameters $\nu$ has also been studied in Section~\ref{sec:accuracy}. Additional information about the Mat\'{e}rn class can be found in~\cite{rasmussen2006gaussian}.

We discuss the computational costs involving matrix vector products $\prior x$ and $\prior^{-1}x$. For stationary or translational invariant covariance kernels with points located on a regular equispaced grid, the computational cost for $\prior x$ can be reduced from $\bigO(m^2)$ using the naive approach, to $\bigO (m\log m)$ by exploiting the connection between Toeplitz structure in 1D or Block-Toeplitz structure in 2D etc, and the Fast Fourier Transform (FFT)~\cite{nowak2003efficient}. For irregular grids, it can be shown that the cost for approximate matrix-vector products (matvecs) involving the prior covariance matrix $\prior$ can be reduced to $\bigO(m\log m)$ using Hierarchical matrices~\cite{saibaba2012efficient} or $\bigO(m)$ using ${\cal{H}}^2$-matrices or Fast Multipole Method (FMM)~\cite{ambikasaran2013large}. For forming matvecs $\prior^{-1}x$, we use an iterative solver such as GMRES with a preconditioner that employs approximate cardinal functions based on local centers and special points~\cite{beatson1999fast}. The cost of constructing the preconditioner is $\bigO(m)$ or $\bigO(m\log m)$. Assuming that the number of iterations is independent of the size of the system, the cost of forming $\prior^{-1}x$ is also $\bigO(n_\text{iter}m)$ or $\bigO(n_\text{iter}m\log m)$, where $n_\text{iter}$ is the number of iterations required to converge to the desired tolerance. In conclusion, the cost for forming $\prior x$ and $\prior^{-1} x$ is $\bigO(m\log m)$.

In the case that the operator $h(s) = Hs$ is linear, the resulting posterior PDF in Equation~\eqref{eqn:posteriorpdf} is also Gaussian. The MAP point~\eqref{eqn:mapestimate} can be calculated by introducing auxillary variables $\hat{\xi}$ and $\hat\beta$ which satisfy the following system of equations 
\begin{equation}\label{eqn:linearmap}
\begin{pmatrix} H\prior H^T + \noise & HX \\ (HX)^T & 0  \end{pmatrix} \begin{pmatrix} \hat\xi \\ \hat\beta \end{pmatrix} = \begin{pmatrix} y \\ 0 \end{pmatrix}  
\end{equation}
The MAP estimate can be calculated as $\hat{s} = X\hat\beta + \prior H^T\hat\xi$. Further details on the $\hat\xi-\hat\beta$ formulation and the procedure to compute the solution to the linear system of equations in Equation~\eqref{eqn:linearmap} can be found in~\cite{saibaba2012efficient,saibaba2012application}. Furthermore, the posterior covariance matrix is given by the expression 
\begin{equation}
\label{eqn:posteriorlinear}
 \post \define \quad  \begin{pmatrix} F_{ss} & F_{s\beta} \\ F_{s\beta}^T &F_{\beta\beta}\end{pmatrix}\quad = \quad  \begin{pmatrix} \prior^{-1} + H^T\noise^{-1} H  & \prior^{-1}X \\ X^T\prior^{-1} & X^T\prior^{-1}X\end{pmatrix}^{-1}
\end{equation}

\section{Posterior Covariance approximation}\label{sec:posteriorapprox}

In~\cite{flath2011fast,bui2012extreme}, the approach taken to approximate the posterior covariance is to compute the dominant modes of a prior preconditioned data misfit Hessian $\prior^{1/2}H_\text{red}\prior^{1/2}$, which combines the information from the data misfit portion of the Hessian $H$ and the prior covariance matrix $\prior$. Here, $\red \define -\nabla_{ss}  \log p(y|s) = H^T\noise^{-1}H$, where $p(y|s)$ is the likelihood of the measurements conditioned on the data. The matrix $H_\text{red}$ often has a rapidly decaying spectrum for several ill-posed inverse problems~\cite{flath2011fast,bui2012extreme,bui2012analysis1,bui2012analysis2}. Computing the dominant modes of this eigenvalue problem can be accomplished in a cost that is a small multiple of the cost of simulating the forward (or adjoint) problem, and this multiple is independent of the size of the parameter dimension. Once this low-rank approximation is computed, an approximation to the posterior covariance matrix can be computed as a low-rank update to the prior covariance matrix using the Woodbury matrix identity.  This results in a scalable algorithm for estimating uncertainty in large-scale statistical inverse problems. %

\subsection{Avoiding  $\prior^{1/2}$ and $\prior^{-1}$}\label{sec:avoiding}
As mentioned earlier, the key drawback of~\cite{flath2011fast,bui2012extreme} is that the algorithms for approximating posterior covariance matrix involves computing (matvecs with) the square root of the prior covariance matrix. When the prior covariance matrix is represented by a discrete representation of a (possibly fractional) differential operator~\cite{stuart2010inverse}, the square root can be computed using sparse Cholesky representation that scales as $\bigO (m^{3/2})$ in 2D and $\bigO(m^2)$ for  3D problems. 

Several matrix-free techniques exist in the literature for computing matvecs $\prior^{1/2}x$ and $\prior^{-1/2}x$, that are based on polynomial approximation~\cite{chen2011computing,dietrich1995efficient} or rational approximations and contour integrals~\cite{hale2008computing}. However, the convergence of polynomial approximations is only algebraic when the smallest eigenvalue is close to zero. Rational approximations and contour integral based methods  do not suffer from the same problem, however they require solutions of a number of shifted systems totaling $\bigO (\log\kappa_\text{cond}(\prior)) $, where $\kappa_\text{cond}(\cdot)$ is the condition number defined for a matrix $A$ as $\kappa_\text{cond}(A) \define \normtwo{A}\normtwo{A^{-1}}$. Although the number of systems to be solved is often small, even for ill-conditioned problems, solving each system can be expensive in practice. Recent developments also include symmetric square-root factorizations, such as $\mathcal{H}$-Cholesky~\cite{bebendorf2005hierarchical}. However, they are complicated to implement and have low accuracy for a reasonable computational cost. Therefore, it is highly desirable to develop square-root free algorithms for approximating the posterior covariance matrix.

Our approach, instead, is to consider the dominant eigenmodes of the generalized Hermitian eigenvalue problem (henceforth referred to as GHEP) $H_\text{red}x=\lambda \prior^{-1}x$, which has the same spectrum as the matrix prior preconditioned data misfit Hessian $(\prior^{1/2}H_\text{red}\prior^{1/2})$. However, as we shall show in Section~\ref{sec:randomized}  the dominant eigenmodes can be computed without relying on square-roots of the matrix $\prior$. We consider the generalized eigenvalue problem 
\begin{equation}
\label{eqn:posteriorghep}
 H^T\noise^{-1} H u = \lambda \prior^{-1} u 
\end{equation}
Recall that $\red = H^T\noise^{-1}H$. Since both matrices $H_\text{red}$ and $\prior$ are symmetric and $ \prior$ is symmetric positive definite, we have the following eigendecomposition
\begin{equation}
\label{eqn:geneigendecomposition}
H_\text{red} = \prior^{-1} U \Lambda U^T \prior^{-1} \qquad U^T\prior^{-1}U = I
\end{equation}
where $U$ is the matrix of eigenvectors obtained by the solution of Equation~\eqref{eqn:posteriorghep} and $\Lambda$ is a diagonal matrix with eigenvalues. Equation~\eqref{eqn:posteriorghep} does not have square-roots but still has $\prior^{-1}$. We further make the following variable transformation, $y = \prior^{-1} x$. Therefore, Equation~\eqref{eqn:posteriorghep} becomes 

\begin{equation}
\label{eqn:posteriorghep2}
 \prior H^T\noise^{-1} H \prior y = \lambda \prior y 
\end{equation}
Note that this transformation still preserves the same eigenvalues $\lambda$. After solving for the eigenpair $(\lambda,y)$ the eigenvector of Equation~\eqref{eqn:posteriorghep} can be recovered by the following transformation $x = \prior y$. As a result, we have completely removed the need for inverting $\prior$ or forming square roots. We now consider efficient solvers for computing the dominant eigenmodes of the GHEP~\eqref{eqn:posteriorghep}.

\subsection{A randomized algorithm for approximating posterior covariance}\label{sec:randomized}
In Section~\ref{sec:avoiding}, we described an efficient representation of the posterior covariance matrix, based on the dominant eigenmodes of the eigenvalue problem described in Equation~\eqref{eqn:posteriorghep}. Equation~\eqref{eqn:posteriorghep} is a GHEP, it is a special case of a generalized eigenvalue problem $Ax = \lambda Bx$ with A being Hermitian and B being Hermitian positive definite. Several efficient methods exist for solving the GHEP. These include approaches based on power and inverse iteration methods, Lanczos based methods and Jacobi-Davidson method. For a good review on this material, please refer to~\cite[chapter 5]{bai1987templates} and~\cite{saad1992numerical}. A good review of existing software available for solving eigenvalue problems is also available in~\cite{str-6}.

In this paper, we choose to use the square-root free randomized algorithms that we have recently developed in~\cite{saibaba2013randomized}. Randomized algorithms are gaining popularity because 1) they are easy to implement, 2) they are well suited for computationally intensive problems and modern computing environments, and 3) have well-established error bounds and computational costs. The error analysis suggests that the randomized algorithm developed in~\cite{saibaba2013randomized} is most accurate when the generalized singular values of the matrix $B^{-1}A$ decay rapidly. In Section~\ref{sec:avoiding}, we have argued that the generalized eigenvalues of the GHEP~\eqref{eqn:posteriorghep} are rapidly decaying and therefore, we expect the randomized algorithm to be fairly accurate. The randomized algorithms that we employ are also extremely efficient because they avoid forming expensive matvecs with $B^{1/2}$ and $B^{-1/2}$. In addition, they produce a symmetric low-rank representation with B-orthonormal eigenvectors. This symmetry in low-rank representation has been exploited in deriving an efficient, symmetric representaion of the posterior covariance matrix in Section~\ref{sec:approx}.

We briefly review the randomized algorithm described in~\cite{saibaba2013randomized} for computing dominant eigenmodes of the GHEP $Ax = \lambda Bx$. In the context of solving the problem~\eqref{eqn:posteriorghep}, we have \[ A\define \prior H^T\noise^{-1}H \prior  \quad \text{and} \quad B \define \prior\]  The key observation is that the matrix $C \define B^{-1}A = \red\prior$ is symmetric with respect to the $B$-inner product $<x,y>_B = y^TBx$. Suppose we wanted to compute the $k$ largest generalized eigenpairs of $Ax=\lambda B x$. The randomized Algorithm~\ref{alg:randomizedghepuncert} calculates a matrix $Q$, which is $B$-orthonormal and approximately spans the column space of $C$, i.e. satisfies the following error bound $\normB{(I-QQ^TB)C} \leq \varepsilon $. Given such a matrix $Q$, it can be shown that $\normB{ A - (BQ) (Q^TAQ) (BQ)^T} \leq 2\varepsilon$, i.e. $A\approx (BQ) (Q^TAQ) (BQ)^T$. As a result, a symmetric rank-$k$ approximation can be computed, from which the approximate eigendecomposition can be computed.

To produce a symmetric rank-$k$ approximation, the algorithm proceeds as follows: first, we sample a matrix with entries randomly chosen from $\normal (0,1)$, $\Omega \in \mathbb{R}^{n\times r}$. We choose $r = k +p$, where $p$ is an oversampling factor, which we choose to be $20$. Numerical evidence for this choice of oversampling factor has been provided in~\cite{halko2011finding,saibaba2013randomized}. Form $\bar{Y} = A\Omega$ and $Y = B^{-1}AQ$. Then, we compute the QR factorization of  $Y = QR$ such that $Q^TBQ = I$. This can be accomplished by modified Gram-Schmidt algorithm with $\langle\cdot,\cdot\rangle_{B}$ inner products. Then, we form $T\define Q^TAQ$ and compute its eigenvalue decomposition $T = S\Lambda S^T$. We then have the approximate generalized eigendecomposition 
\[ A \approx U\Lambda U^T \qquad U = QS\]
Here, $U$ is also $B$-orthonormal. The cost of a second round of matvecs with $A$ while computing $T$ can be avoided using the following approximation
\[ \Omega^T\bar{Y} = \Omega^TA\Omega \approx (\Omega^T BQ) T (Q^TB\Omega)\]
Therefore, $T$ can be computed as $T \approx (\Omega^T BQ)^{-1}(\Omega^T\bar{Y}) (Q^TB\Omega)^{-1}$. The results are compactly summarized in Algorithm~\ref{alg:randomizedghepuncert}. For an efficient implementation we describe the algorithm using $A = \prior H^T\noise^{-1}H \prior$ and $B = \prior$ and $B^{-1}A = H^T\noise^{-1}H \prior$.

\begin{algorithm}[!ht]
\begin{algorithmic}[1]
 \REQUIRE  matrices $\red, \prior  \in \mathbb{R}^{n\times n}$ and $\Omega \in \mathbb{R}^{n\times (k+p)}$ is a Gaussian random matrix. Here, k is the desired rank, $p\sim 20$ is an oversampling factor. 
\STATE  Compute $ Y =   \red \prior \Omega$
\STATE  Compute QR - factorization of $Y = QR$ s.t. $Q^T\prior Q = I$
\STATE  Form $ T = Q^T\prior\red\prior Q$ or $(\Omega^T BQ)^{-1}(\Omega^T\bar{Y}) (Q^TB\Omega)^{-1}$ 
\STATE  Compute the eigenvalue decomposition $T = S\Lambda S^T$. Keep the $k$ largest eigenmodes.
\STATE  \textbf{Return:} Matrices $U \in \mathbb{R}^{n\times k}$ and $\Lambda\in \mathbb{R}^{k\times k}$ 
 that satify
\[ \red  \approx  (\prior^{-1}U) \Lambda (\prior^{-1}U)^*\qquad\text{with}\qquad U = \prior QS \]
\end{algorithmic}
\caption{Randomized algorithm for computing dominant eigenmodes  of $\red x = \lambda \prior^{-1}x$ }
\label{alg:randomizedghepuncert}
\end{algorithm}

The efficiency and accuracy of this algorithm has been studied in detail in~\cite{saibaba2013randomized}. Here we summarize the main conclusions. The error in the low-rank approximation is $$\normB{(I-QQ^TB)B^{-1}A} \leq c\normtwo{B^{-1}}\sigma_{B,k+1}(B^{-1}A)$$
where $c$ is a constant that depends on $n,k$ and $p$ and is independent of the spectrum of the matrices, $\sigma_{B,k+1}$ is the $(k+1)$-th generalized singular value of the matrix $B^{-1}A$. Since the generalized singular values are not known in advance, a randomized estimator for the error in the low-rank representation is also proposed and analysed in~\cite{saibaba2013randomized}. Given error in the low-rank representation $\normB{(I-QQ^TB)B^{-1}A}$, it can be shown that the error in approximating the true eigenvalue and eigenvector satisfies the following error bounds  
\[|\lambda-\tilde{\lambda}| \leq \min\{2\varepsilon,\frac{4\varepsilon^2}{\delta}\}\qquad \sin\angle_B (u,\tilde{u}) \leq \frac{2\varepsilon}{\delta}\]

where $\delta = \min_{\lambda_i\neq\lambda} |\tilde{\lambda}-\lambda_i|$ is the gap between the approximate eigenvalue $\tilde{\lambda}$ and any other eigenvalue and $\angle_B(x,y) = \arccos \frac{|<x,y>_B|}{\normB{x}\normB{y}}$. This result states that the accuracy in the eigenvalue/eigenvector calculations depends not only on the accuracy of the low-rank representations but also on the spectral gap $\delta$. When the eigenvalues are clustered, the spectral gap is small and the eigenvalue calculations are accurate as long as the error in the low-rank representation is small. However, in this case the resulting eigenvector calculations maybe inaccurate because the parameter $\delta$ appears in the denominator for the approximation of the angle between the true and approximate eigenvector. This result has consequences in the accuracy approximation of the posterior covariance and will be discussed in Section~\ref{sec:accuracy}.

\subsection{Posterior covariance approximation}\label{sec:approx}
Let us define the matrix $F_{ss} \define \left(H^T\noise^{-1}H + \prior^{-1}\right)^{-1}$. We follow the approach in~\cite{flath2011fast,martin2012stochastic} to derive an approximation to $F_{ss}^{-1}$ as a low-rank update to the prior covariance matrix $\prior$. Plugging in the approximate eigendecomposition for $H_{red}$ from Equation~\eqref{eqn:geneigendecomposition}, we obtain 
\[  F_{ss} \quad =  \quad  \left(\prior^{-1} U\Lambda U^T\prior^{-1}U  + \prior^{-1}\right)^{-1}\]
Using the Sherman-Morrison-Woodbury formula, we can derive the following expression
\begin{align}\label{eqn:fssinvapprox}
F_{ss}^{-1}  \quad =  & \quad \prior -  \prior \prior^{-1} U \left(\Lambda^{-1}  + U^T\prior^{-1}U\right)^{-1}U^T \prior^{-1} \prior \\ \nonumber 
                =  & \quad  \prior - U (\Lambda^{-1} + I)^{-1}U^T   \\ \nonumber 
		= & \quad  \prior -  U_k D_k U_k^T  + \bigO \left(\frac{\lambda_{k+1}}{1+\lambda_{k+1}}\right)   \nonumber
\end{align}
where $D_k \define \text{diag} ( \frac{\lambda_i}{1 + \lambda_i} ) \in \mathbb{R}^{k\times k}$ for $i=1,\dots,k$. The error in the approximation to the matrix $F_{ss}$ and its inverse can be established by the following inequalities

\begin{align}
\label{eqn:inverseerror}
\norm{F_{ss} - \left(\prior^{-1} + \prior^{-1}U_k\Lambda_k U_k^T\prior^{-1} \right)}{M} \quad & \leq \quad  \lambda_{k+1} \normtwo{\prior^{-1}}\\ \nonumber 
\norm{F_{ss}^{-1} - \left(\prior -   U_k D_k U_k^T \right)}{M} \quad & \leq \quad\frac{\lambda_{k+1}}{1+\lambda_{k+1}} \normtwo{\prior}    \nonumber  
\end{align}

The inequality holds for both the cases that $M =\prior$ and $M =\prior^{-1}$. Here, the induced matrix norm is defined as $\norm{A}{M} \define \max_{\norm{x}{M} = 1} \norm{Ax}{M}$, defined for symmetric positive definite matrices $M$.

For many ill-posed inverse problems, the approximate numerical rank $k$ of the prior preconditioned Hessian $\prior^{1/2}\red\prior^{1/2} $ is small and independent of the problem size, i.e. the number of unknowns. Since this matrix has the same eigenvalues as the GHEP $\red x = \lambda \prior^{-1}x$, we are justified in retaining only the eigenmodes corresponding to the largest eigenvalues. We can exploit this low-rank structure to develop scalable algorithms to approximate the posterior covariance matrix. The generalized eigendecomposition combines information from the prior and the reduced Hessian and takes advantage of the eigenvalue decay in one (or both) matrices - when the reduced Hessian has rapidly decaying eigenvalues, or the prior is of smoothing type.  Analytical evidence for the eigenvalue decay of the reduced Hessian $H_\text{red}$ is provided in~\cite{flath2011fast} in the context of advection-diffusion based inverse problems and in~\cite{bui2012analysis1,bui2012analysis2} for inverse scattering problems. For the case of prior covariance matrices, the eigenvalue spectrum is known to decay rapidly when the covariance kernels are smooth~\cite{schwab2006karhunen}. The $k$ retained eigenvectors are the modes along which the parameter field is informed by a combination of the data and the prior. Typically the data and prior are informative about the low-frequency modes and as a result local information and fine scale information is hard to recover from the data. The inequalities in~\eqref{eqn:inverseerror} suggests choosing $k$ such that  $\lambda_{k+1} \ll 1$.

Finally, the posterior covariance matrix can be constructed by using the following decomposition in Equation~\ref{eqn:posteriorlinear}
\begin{align} \nonumber
\post \quad = & \quad \begin{pmatrix} I & -F_{ss}^{-1}F_{s\beta} \\ & I \end{pmatrix} \begin{pmatrix} F_{ss}^{-1}  & \\ & (F_{\beta\beta}^{-1} - F_{s\beta}^TF_{ss}^{-1} F_{s\beta})^{-1} \end{pmatrix} \begin{pmatrix} I & \\ - F_{s\beta}^TF_{ss}^{-1} & I \end{pmatrix}  \\ \label{eqn:fisherinverse}
= & \quad \begin{pmatrix} F_{ss}^{-1}  +GS^{-1} G^T  &  -GS^{-1} \\ -S^{-1} G^T  & S^{-1}\end{pmatrix}  \end{align}
where $S\define (F_{\beta\beta}^{-1} - F_{s\beta}^TF_{ss}^{-1} F_{s\beta})$ and $G \define F_{ss}^{-1}F_{s\beta}$. 

An approximation to the posterior covariance matrix~\eqref{eqn:fisherinverse} can be constructed by plugging in the approximation to $F_{ss}^{-1}$ which has been derived in Equation~\eqref{eqn:fssinvapprox}. It should be noted that computing and storing the posterior covariance matrix or its inverse is infeasible for very large problem sizes. However, neither computing nor storing the matrix $\post$ is necessary, nor recommended.  Computing $F_{ss}^{-1}x$ can be performed as
\[ F_{ss}^{-1}x \quad \approx \quad  \prior x - U_kD_kU_k^T x \qquad \text{Cost :} \quad \bigO (m\log m + km)\]
By the same argument, forming $G$ and $S$ is approximately the same cost as forming matrix-vector product with $F_{ss}^{-1}$. Having computed $G$ and $S$, the matrix-vector product involving $\post$ is dominated by the cost of computing $F_{ss}^{-1}x$ and can be computed efficiently in $\bigO (m\log^\gamma m + rm)$.

\section{Application: Ray Tomography}\label{sec:raytomog}
\subsection{Setup}
In this application, we consider a synthetic  setup of cross-well tomography~\cite{ambikasaran2013large}. The goal is to the image the slowness in the medium where slowness is defined as the reciprocal of seismic velocity. As a first order approximation, the seismic wave is modeled as traveling along a straight line from the sources to the receivers without reflections or refractions. Each source-receiver pair generates one measurement and therefore, there are $n = n_\text{rec} n_\text{sou}$ measurements. Here $n_\text{rec}$ is the number of receivers and $n_\text{sou}$ are the number of sources. In this application, we pick $n_\text{sou} = 20$ and $n_\text{rec} = 50$. The domain is discretized into $m$ cells $\sqrt{m}\times\sqrt{m}$ and within each cell, the slowness is assumed to be constant. Therefore, the time taken from the source to the receiver is a weighted sum of the slowness in the cell, weighted by the length of the ray within the cell. The inverse problem can be stated as follows: given measurements of time delay between each source-receiver pair, what is the slowness of the medium? We assume that the travel times are corrupted by Gaussian noise, so that the measurement takes the following form 

\begin{equation} \label{eqn:lineartomography}
y = Hs + v \qquad v \sim \normal (0,\noise )
\end{equation}
where $y$ are the observed (synthetic) travel times, $s$ is the slowness that we are interested in imaging and $H$ is the measurement operator, whose rows correspond to each source-receiver pair and are constructed such that their inner product with the slowness would result in the travel time. Constructed as above results in $H$ as a sparse matrix with $\bigO (n\sqrt{m})$ non-zero entries - each row has $\bigO (\sqrt{m})$ entries and there are $n$ measurements. As a ``true field'', we take a realization from a Gaussian random field $s_\text{true} \sim \normal (X\beta_\text{true},\prior)$, with $X = [1,\dots,1]^T$, $\beta_\text{true} = 5\times 10^{-3}[s/m]$ and covariance matrix $\prior$ is constructed according to the covariance kernel 
\begin{equation}\label{eqn:raycovariance}
\kappa(\bx,\textbf{y}) = \theta\exp\left( -\frac{\normtwo{\bx-\textbf{y}}}{L} \right)
\end{equation}  
with $\theta = 1\times 10^{-3}[s/m]$ and $L = 100 [m]$. The realization of the stationary Gaussian process is generated using the FFT-based approach in~\cite{dietrich1997fast} by embedding the Toeplitz matrix in a circulant matrix. The matvecs with $\prior$ are also performed with FFT based methods in $\bigO(m\log m)$.

\begin{figure}[!ht]
\centering
\includegraphics[scale = 0.5]{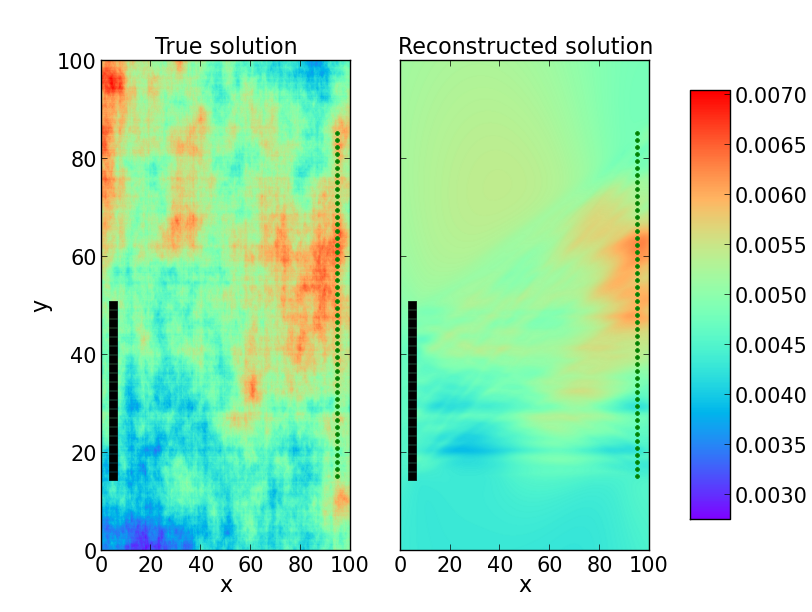}
\caption{(left) True field which is a realization of a Gaussian process with mean $5\times 10^{-3}$[s/m] and covariance kernel~\eqref{eqn:raycovariance}, and (right) Reconstruction using the geostatistical approach on a grid with $1024\times 1024$ points. Black solids correspond to source locations (20) and green dots correspond to the receiver locations (50). The reconstruction error in relative $L^2$ sense was $0.11$. }
\label{fig:tomogrecon}
\end{figure}

\begin{figure}[!ht]
\centering
\includegraphics[scale = 0.5]{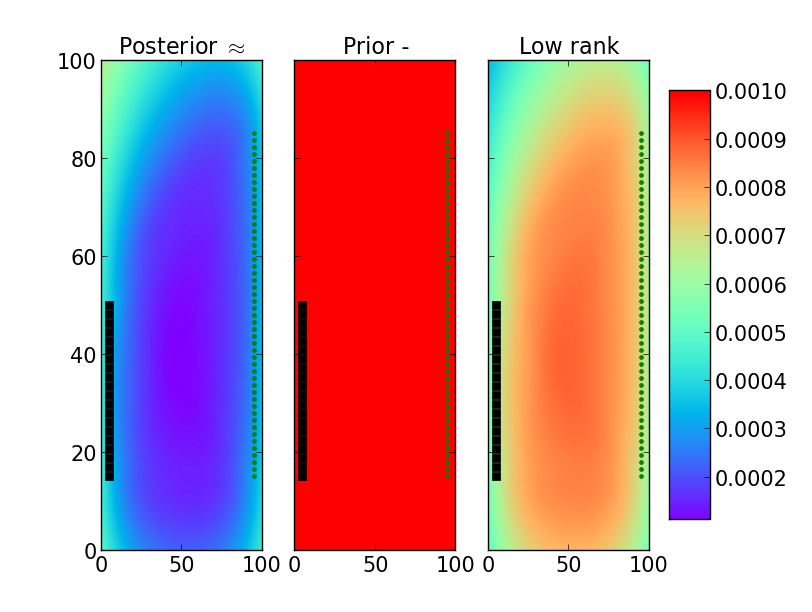}
\caption{Visualization of the diagonals of the posterior covariance matrix. (left) diagonals of the posterior covariance matrix, (middle) diagonals of the prior covariance matrix and (right) the low-rank correction from the $(1,1)$ block in Equation~\eqref{eqn:fisherinverse}. The uncertainty in the estimation of $\beta$ measured as $\trace (S^{-1})$ (see Equation~\eqref{eqn:posteriorlinear}) is $6.32\times 10^{-4}$ [s/m].}
\label{fig:tomogdiags}
\end{figure}

\begin{figure}[!ht]
\centering
\includegraphics[scale = 0.5]{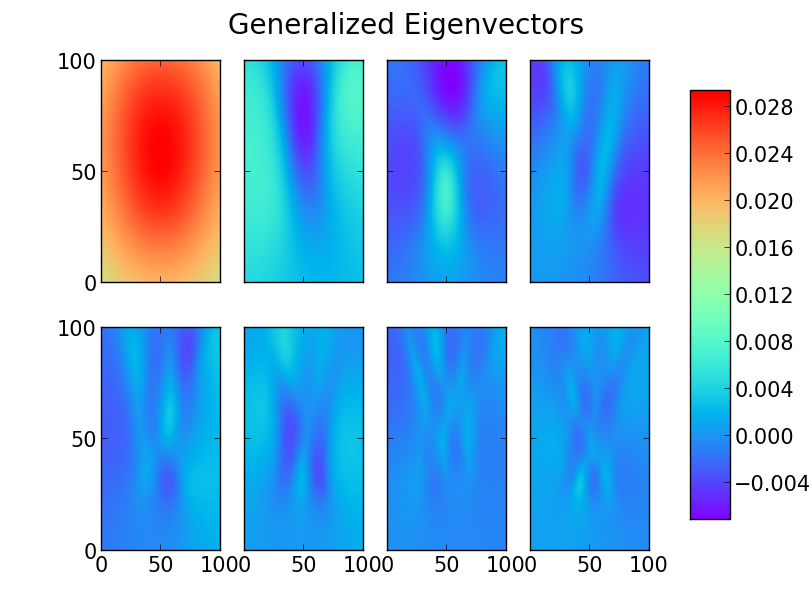}
\caption{Eigenvectors $1,3,5, 7,10,12,18$ and $20$ of the generalized eigenvalue problem~\eqref{eqn:geneigendecomposition} solved using Algorithm~\ref{alg:randomizedghepuncert}. }
\label{fig:rayeigenvecs}
\end{figure}

\begin{figure}[!ht]
\centering
\includegraphics[scale = 0.5]{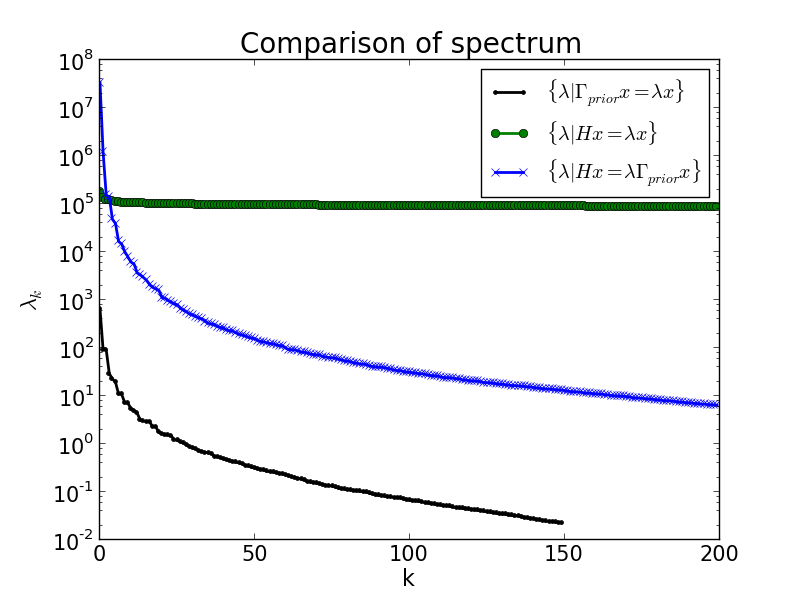}
\caption{Comparison of the eigenvalues of the data misfit Hessian $H_\text{red} = H^T\noise^{-1}H$, the prior covariance matrix $\prior$ and the GHEP in the Equation~\eqref{eqn:posteriorghep}. As can be seen clearly, the eigenvalues of the GHEP that combines information from both the prior $\prior$ and the data misfit Hessian $\red$ decay rapidly. The grid size used in this calculation was $1024\times 1024$. }
\label{fig:tomogeigencompare}
\end{figure}

To simulate experimental error, we add Gaussian noise of 0.1$\%$ to the measurements. The system of Equations~\eqref{eqn:linearmap} is solved using an iterative solver, which is chosen to be restarted GMRES (50), and the matvecs can be formed in a matrix-free manner~\cite{saibaba2012efficient,saibaba2012application}. As a preconditioner, we use the method developed in~\cite[Section 3.3]{saibaba2012efficient}. The preconditioner requires computing a low-rank decomposition of $\prior$ which can be computed using any eigenvalue solver. We use the randomized GHEP algorithm described in~\ref{sec:randomized} with $A = \prior$ and $B=I$, with $\sim 150$ eigenpairs used to form the low-rank decomposition and an oversampling factor $p=20$. The other steps in the construction of the preconditioner are identical. 
To simulate experimental error, we add Gaussian noise of 0.1$\%$ to the measurements.  The iterative solver converged in $138$ iterations. The $L^2$ error in the reconstruction was $11.08\%$.  Figure~\ref{fig:tomogrecon} shows the true field which is a realization of a Gaussian process with mean $5\times 10^{-3}$[s/m] and covariance kernel~\eqref{eqn:raycovariance} and the reconstructed field obtained by solving Equation~\eqref{eqn:linearmap}. The grid size was $1024\times 1024$  and the relative reconstruction error in the l$_2$ sense was $0.11$. Figure~\ref{fig:tomogdiags} visualizes the posterior variance obtained as the diagonals of the posterior covariance matrix calculated using the formula $\post = \prior -U_kD_kU_k^T$. The error in the uncertainty associated with estimation of $\hat\beta$ measured as $\trace (S^{-1})$ (see Equation~\eqref{eqn:posteriorlinear}) is $6.32\times 10^{-4}$ [s/m]. The visualization a few selected eigenvectors of the generalized eigenvalue problem~\eqref{eqn:geneigendecomposition} solved using Algorithm~\ref{alg:randomizedghepuncert} can be seen in Figure~\ref{fig:rayeigenvecs}. 

\subsection{Variance approximation}
\label{sec:accuracy}
In this section, we study the effect of truncation of the spectrum and the choice of prior on the accuracy of the variance of the reconstruction, that is the diagonals of the posterior covariance matrix. We consider three different covariance kernels chosen from the Mat\'{e}rn class of covariance kernels corresponding to parameters $\nu = 1/2,3/2,5/2$, i.e. 

\begin{equation}\label{eqn:maternnu}
\kappa_\nu (\bx,\textbf{y}) = \left\{ \begin{array}{ll} \theta\exp(-r/L) &  \quad \nu = 1/2 \\ \theta(1 + \sqrt{3}r/L)\exp(-\sqrt{3}r/L) & \quad \nu =3/2  \\  \theta (1 + \sqrt{5}r/L + 5(r/L)^2)\exp(-\sqrt{5}r/L) & \quad \nu =5/2  \end{array}\right.
\end{equation}

\begin{figure}[!ht]
\centering
\includegraphics[scale = 0.35]{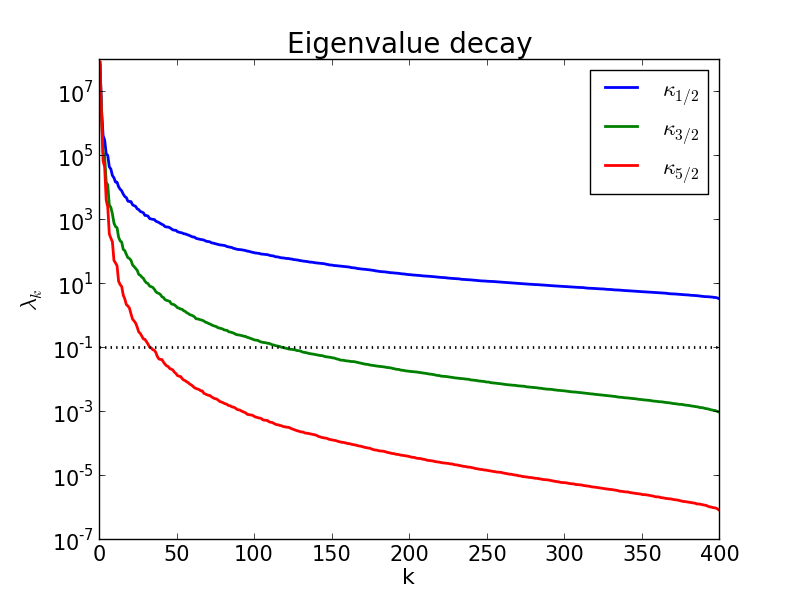}
\includegraphics[scale = 0.35]{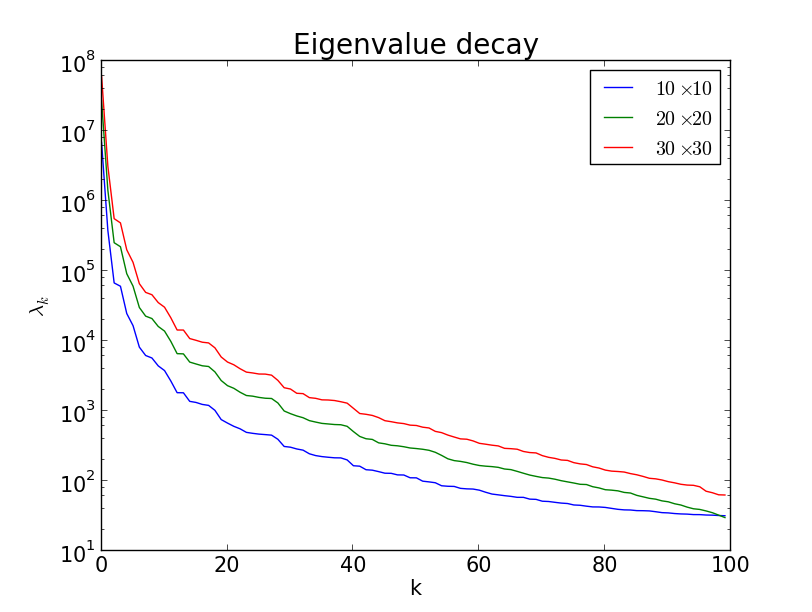}
\caption{(left) The decay of the eigenvalues for different covariance kernels. The unknowns are discretized on a $200\times 200$ grid and the number of measurements are $n_m = 400$. Three different covariance kernels are considered described in Equation~\eqref{eqn:maternnu} with $\theta = 10^{-3}$. (right) The decay of the eigenvalues for different setups with differing number of measurements ranging from $10\times 10$ to $30\times 30$. The covariance kernel was chosen to be $\kappa_{1/2}$ described in Equation~\eqref{eqn:maternnu} and the unknowns were discretized on a $200\times 200$ grid. Only the top $100$ eigenvalues are displayed. }
\label{fig:eigdecaycovariance}
\end{figure}

The comparison of the rate of decay of the eigenvalues of the matrices $\prior$, $\post$ and the GHEP $H_\text{red}x = \lambda \prior^{-1}x$ is illustrated in Figure~\ref{fig:tomogeigencompare}. What does affect the spectrum of the generalized eigenvalue problem for a given measurement setup is the choice of the covariance kernel used in the reconstruction. Here we illustrate this by choosing a ray tomography setup with $n_m=400$ measurements and unknowns discretized on a $200\times 200$ grid. For the Mat\'{e}rn class of covariance kernels, by increasing the parameter $\nu$, the underlying stochastic process becomes more smooth and the rate of decay of the eigenvalues increases~\cite{rasmussen2006gaussian,schwab2006karhunen}. This affects the number of eigenvalues retained to maintain the accuracy of the posterior covariance matrix - smoother the kernel, fewer eigenvalues need to be retained. The number of retained eigenvalues is $\min\{n_m,k(\varepsilon)\}$, where $k(\varepsilon)$ is chosen such that the retained eigenvalues satisfy the cutoff $\lambda > \varepsilon$. A typical choice is $\varepsilon = 0.1$.  This is illustrated in Figure~\ref{fig:eigdecaycovariance}. We see that for the covariance kernel corresponding to $\nu=1/2$, the eigenvalue decay is not rapid enough to satisfy the cutoff and as a result, all the non-zero eigenvalues that correspond to the number of measurements $n_m = 400$. Physical parameters also control the decay of eigenvalues and as a result, the number of eigenmodes that are retained. Setups that have more number of measurements due to increased number of sources and receivers, have a slower rate of decay of eigenvalues because more information is propagated through the rows of the measurement operator. This is also illustrated in Figure~\ref{fig:eigdecaycovariance}.

\subsubsection{Accuracy of variance approximation}	
\begin{figure}[!ht]
\centering
\includegraphics[scale = 0.35]{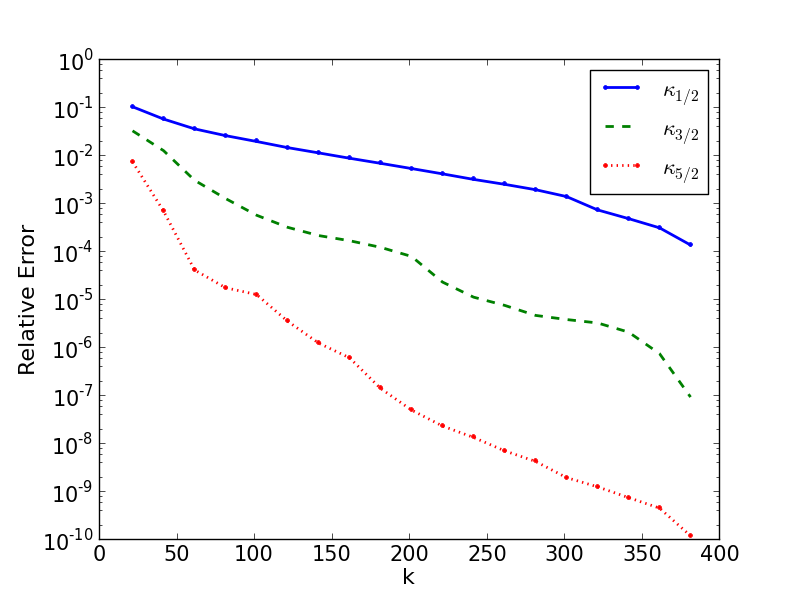}
\caption{Dependence of the relative error in the approximation to the variance computed as a function of number of retained eigenvalues retained in the approximation. The relative error computed with the exact variance, where relative error is $\sum_k |\text{Var}_k^\text{exact} - \text{Var}_k^\text{approx}| / \sum_k |\text{Var}_k^\text{exact}|   $ }
\label{fig:eigerror}
\end{figure}

In Figure~\ref{fig:eigerror}, we plot the relative error in the computation of the approximate variance as a function of the number of eigenvalues retained for different prior covariance kernels. The error  is defined to be the relative error computed with the exact variance.  We observe that the error decreases with increasing number of eigenvalues retained. Consequently, inclusion of a larger number of eigenvalues presents a situation of diminishing returns in terms of the accuracy of variance calculations. The result of the computation of the accuracy of the variance with increasing eigenvalues is plotted in Figure~\ref{fig:eigerror}. As can be seen, the error in the approximation to the variance decreases rapidly with increasing number of eigenvalues retained. Furthermore, the accuracy improves dramatically when  considering kernels with higher values of $\nu$,  since eigenvalues of covariance matrices  arising from the Mat\'ern family corresponding to larger $\nu$ decay more rapidly. In our computations, we assumed that the number of measurements are only $400$, this makes the task of computing the exact variance computationally feasible using the method described in~\cite{ambikasaran2013large}. For larger number of measurements, we cannot compute the variance easily and we must resort to the methods described in this paper.  

\subsubsection{Cost of computing variance}
Two factors are necessary to ensure that the cost of computing the variance scales linearly with the dimensionality of the unknowns parameters - the number of eigenvalues retained is independent of the mesh size and the various costs required to capture the dominant eigenmodes scales linearly with the number of unknowns.

\begin{figure}[!ht]
\centering
\includegraphics[scale = 0.35]{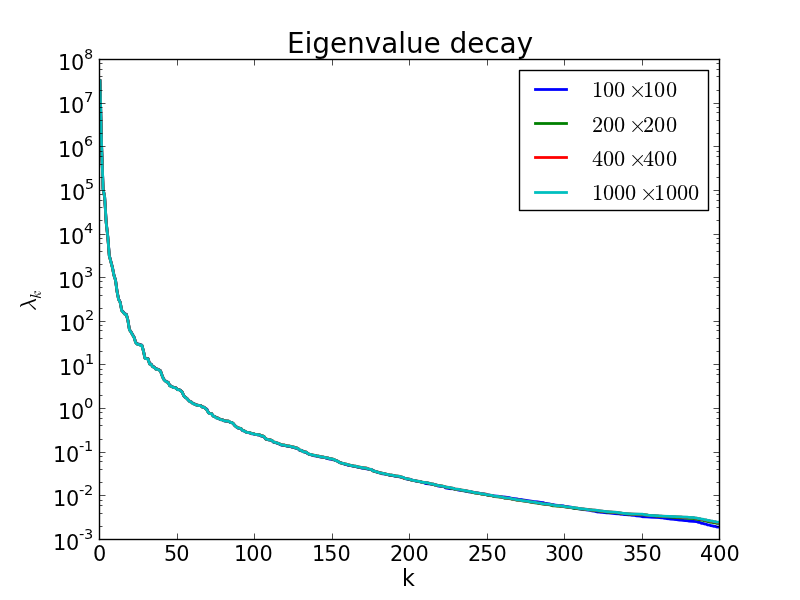}
\includegraphics[scale = 0.35]{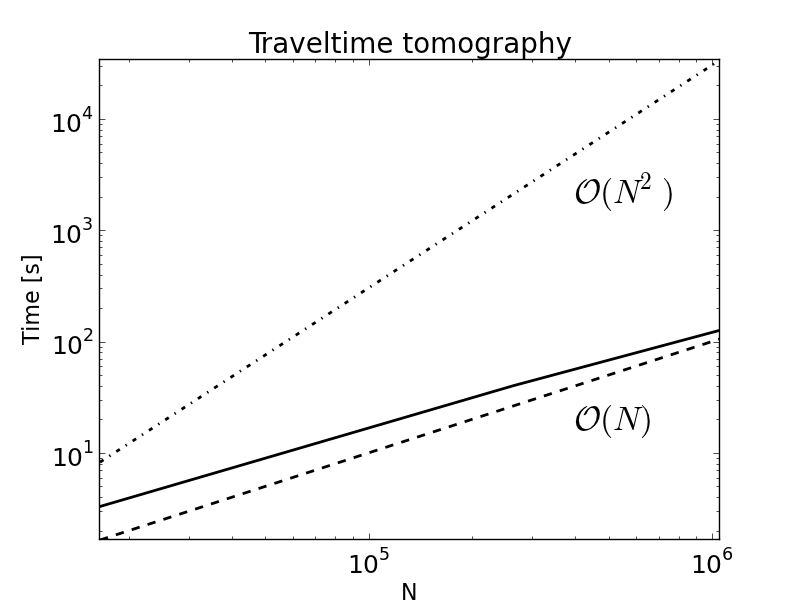}
\caption{(left) The decay of the eigenvalues for the covariance kernel $\kappa(r) = \theta(1 + \sqrt{3}r/L)\exp(-\sqrt{3}r/L)$ where the unknowns are discretized on a number of grids ranging from $100\times 100$ to $1000\times 1000$. (right) The CPU time (solid line) for computing the top $300$ eigenmodes of the ray tomography example with unknowns discretized on grids ranging from $128\times 128$ to $1024\times 1024$ with $n_m = 1000$ measurements. The slopes of the dashed and dotted lines correspond to computational costs of $\bigO(N)$ and $\bigO(N^2)$ respectively. As can be seen, the computational cost scales linearly.}
\label{fig:eigrefinement}
\end{figure}

We address the first issue, that is, the number of eigenvalues retained is independent of the mesh size. Figure~\ref{fig:eigrefinement} shows the dependence of the spectrum of the generalized eigenvalue problem $H_\text{red}x = \lambda \prior^{-1} x$ with the dimension of the input parameter space. We use the same setup just described above, with $400$ measurements. The unknown parameters, in this case the slowness of the medium, is discretized on a number of grids ranging from $100\times 100$ to $1000\times 1000$. The prior is chosen as a member of the Mat\'{e}rn covariance class with $\nu = 3/2$ i.e. $\kappa_{3/2}$, see Equation~\eqref{eqn:maternnu}. As can be seen from the figure, the number of unknowns does not affect the spectrum of the generalized eigenvalue problem and thus, does not affect the number of retained eigenvalues required to approximate the posterior covariance matrix. We address the second issue - computational costs of computing eigendecomposition. The major component in this calculation is the matvecs involving $H_\text{red}$, $\prior $ and $\prior^{-1}$. In Section~\ref{sec:posteriorapprox}, we have already argued that the costs involving the aforementioned matvecs scale linearly or almost linearly with the number of unknowns. Therefore, we expect the cost of the eigendecomposition, and as a result the variance calculations, to scale similarly with the number of unknowns. This is illustrated in Figure~\ref{fig:eigrefinement}. As can be seen from the figure, the computational costs scale linearly with the number of unknowns.

\section{Application: Steady-state Hydraulic Tomography}~\label{sec:ssht}

Hydraulic Tomography (HT) is a technique for estimating the subsurface parameters such as transmissivity by a series of pumping tests, in which pressure (head) is measured and this data is used to reconstruct the parameters of interest by using suitable inversion algorithms. The governing equations are given by the ground water flow equations, which in 2D take the following form

\begin{align}\label{eqn:darcy}
-\nabla (K(x) \nabla \phi)  \quad = & \quad  Q(x - x_s)  \quad&   x \in \Omega  \\ 
\phi  \quad = & \quad 0 \quad  &x \in \partial \Omega\nonumber 
\end{align}
where $x_s$ are the pumping locations, and homogenous Dirichlet boundary conditions are applied everywhere on the boundary $\partial\Omega$.

\begin{figure}[!ht]
\centering
\includegraphics[scale = 0.5]{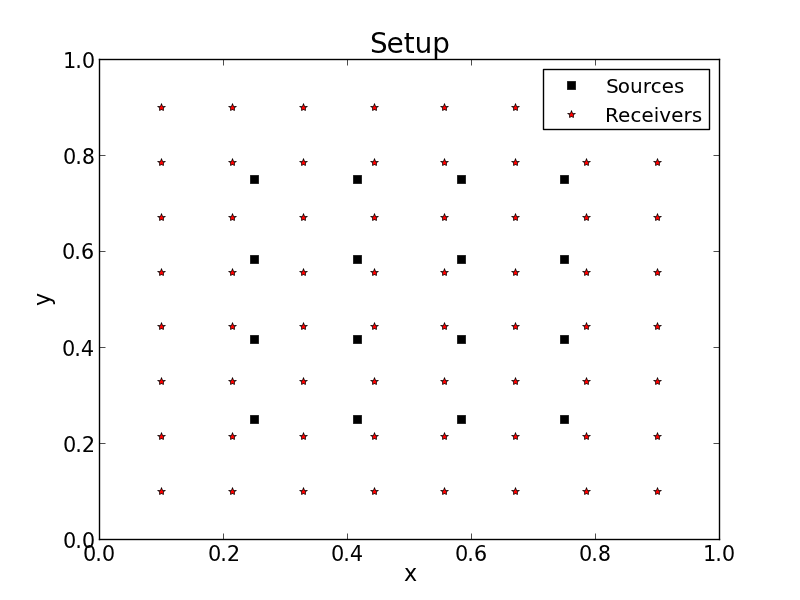}
\caption{Setup of the 16 pumping sources (black solids) and 64 measurement locations (red asterisks). In all, there were $16\times 64 = 1024$ measurements. }
\label{fig:htsetup}
\end{figure}

The inverse problem is to reconstruct the transmissivity field $K(x)$ from discrete measurements of $\phi$ obtained by repeated pumping tests at several pumping locations. In order for the system of Equations~\eqref{eqn:darcy} to be well posed, the reconstructions of $K(x)$ need to be positive and therefore, we make a log-transformation $s(x) \define \log K(x)$ to ensure positivity. This transformation makes the problem of reconstruction a quasi-linear inverse problem. The setup of the pumping tests is provided in Figure~\ref{fig:htsetup}; the black squares indicate the locations at which water is pumped and the head response is computed at the receiver locations (indicated by red asterisks). 

The reconstruction of  transmissivity field $K(x)$ from discrete measurements of $\phi$ is a nonlinear ill-posed inverse problem. For nonlinear inverse problems, unlike linear inverse problems, the posterior distribution is no longer Gaussian. The measurement operator is linearized about the current estimate $s_k$ to obtain 
\[ h(s) \approx h(s_k) + J(s- s_k)  \qquad J = \at{\frac{\partial h(s)}{\partial s}}{s=s_k}\]
The MAP point is the minimizer of the negative log-likelihood posterior PDF described in Equation~\eqref{eqn:mapestimate} and can be computed using the quasi-linear geostatistical approach~\cite{kitanidis1995quasi}. Essentially this involves solving a sequence of regularized linear least squares problems until the algorithm converges.  The system of Equations~\eqref{eqn:inversion} is solved  using restarted GMRES. The prior covariance matrix $\prior$ and the Jacobian  $J_k$ are not constructed explicitly. Matvecs involving $\prior$ are handled in $\bigO(m\log m)$ using the FFT based approach or using $\mathcal{H}$-matrices, whereas the matvecs involving the Jacobian are handled using the approach described in~\cite{haber2000optimization}. In order to accelerate the convergence of GMRES, we use a preconditioner developed in~\cite{saibaba2012efficient}. The algorithm for computing the MAP estimate is summarized in Algorithm~\ref{alg:quasi}.

\begin{algorithm}[!ht]
\caption{Quasi-linear Geostatistical Approach}
\label{alg:quasi}
\begin{algorithmic}[1]
\WHILE {Not converged}
\STATE  Solve the system of equations, 
\begin{equation} \label{eqn:inversion}
\left( \begin{array}{cc} 
        J_k \prior J_k^T  + \noise & J_kX \\	
	 \left(J_k X\right)^T & 0 
       \end{array}
       \right) 
       \left( \begin{array}{c}{\xi_{k+1}}\\ {\beta_{k+1}} \end{array} \right) = 
\left( \begin{array}{c} y - h({s}_k) + J_k{s}_k \\ 0 \end{array} \right) 
\end{equation}

where the Jacobian is defined as $J_k = \at{\frac{\partial h(s)}{\partial s}}{s=s_k}$.

\STATE The update $s_{k+1}$ is computed by, 
\begin{equation} 
 s_{k+1} = X \beta_{k+1} + \prior J_k^T \xi_{k+1} 
\end{equation}

\COMMENT {If necessary, add a line search}
\ENDWHILE

\end{algorithmic}
\end{algorithm}

Once we have the MAP estimate, we approximate the measurement operator by a linearizing about the MAP point so as to approximate the posterior distribution about the MAP estimate by a Gaussian distribution. This is a reasonable approximation when the measurement operator is nearly linear. This is a useful approximation even when the posterior distribution is highly nonlinear, since the Gaussian approximation is frequently used as a proposal distribution for MCMC algorithms used to explore the posterior~\cite{petra2013computational}. The approximate posterior distribution $\tilde{p}(s,\beta| y)$ can be expressed as 
\begin{equation}\label{eqn:postcovquasi}
 \tilde{p}(s,\beta| y) \sim \normal \left( \begin{bmatrix} \hat{s} \\ \hat\beta\end{bmatrix} , \post \right) 
\qquad \post =\begin{pmatrix} \prior^{-1} + J^T\noise^{-1} J  & \prior^{-1}X \\ (\prior^{-1}X)^T & X^T\prior^{-1}X\end{pmatrix}^{-1}
\end{equation}
 
where $\hat{s}$ and $\beta$ are the parameters evaluated at the MAP point, and $J$ is the Jacobian calculated at the MAP point $\hat{s}$.

Computing the  posterior covariance matrix $\post$ or its inverse explicitly is computationally infeasible for large scale problems. This is because (1) constructing the Jacobian explicitly requires the solution of several (possibly time-dependent) forward or adjoint partial differential equations (2) the covariance matrix is typically dense, and therefore computing matrix products and inverses involving the covariance matrix is nearly impossible in terms of storage and computational cost and finally, (3) the resulting matrix is dense and its storage can be infeasible for large problems. The approximation of the posterior covariance matrix is approximated using the same randomized approach described in Section~\ref{sec:posteriorapprox}.

We discuss the implementation of these ideas to a concrete application of steady-state hydraulic tomography. For the reconstruction, we use a covariance kernel, that is part of the Mat\`ern family, and is given by 
\begin{equation}\label{eqn:htcovariance}
\kappa(x,y) = \exp\left( -\frac{\normtwo{x-y}}{4L} \right)
\end{equation}

\begin{figure}[!ht]
\centering
\includegraphics[scale = 0.7]{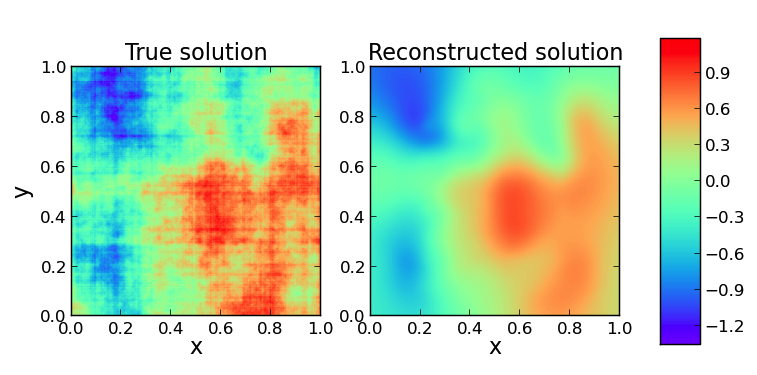}
\caption{(left) True field which is a realization of a Gaussian process with zero mean and covariance kernel~\eqref{eqn:htcovariance}, and (right) reconstruction using the geostatistical approach on a grid with $401\times 401$ points. The reconstruction error in relative $L^2$ sense was $0.29$. }
\label{fig:htrecon}
\end{figure}

%
\begin{figure}[!ht]
\centering
\includegraphics[scale = 0.5]{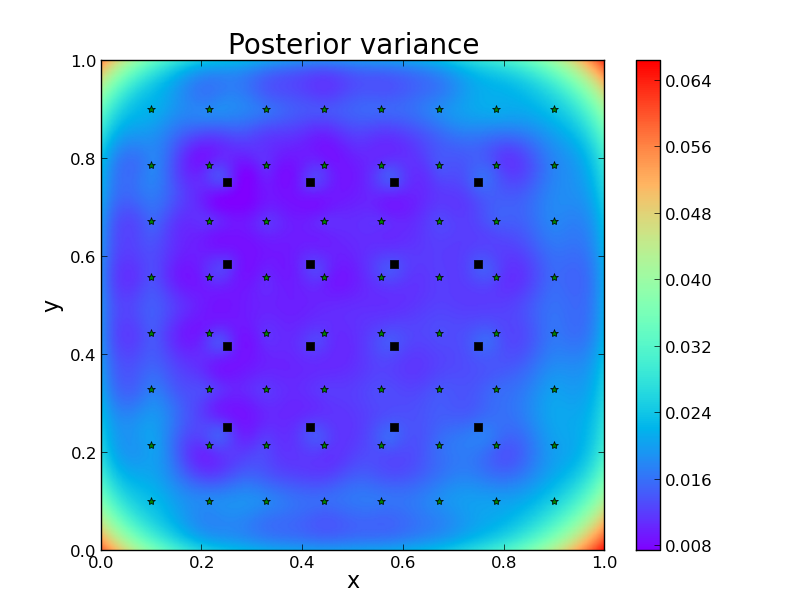}
\caption{Visualization of the diagonals of the posterior covariance matrix in Equation~\eqref{eqn:postcovquasi}. The computation was calculated with the Jacobian evaluated at the MAP estimate.}
\label{fig:htpost}.
\end{figure}

To simulate experimental error, we add Gaussian noise of 0.1$\%$ to the measurements. The MAP estimate is computed by solving the optimization problem~\eqref{eqn:mapestimate} using the Gauss-Newton algorithm that was described in Algorithm~\ref{alg:quasi}. At each iteration the system of equations is solved using restarted GMRES (50) using the same preconditioner described in~\cite{saibaba2012efficient}. The low-rank decomposition is computed in the same fashion described in Section~\ref{sec:raytomog} with $30$ terms in the low-rank approximation. The low-rank decomposition is computed only once; however, the preconditioner is re-built every Gauss-Newton iteration. The iterative solver took about $16$ iterations on average to converge to a relative tolerance of $10^{-7}$. The Gauss-Newton solver took about $6$ iterations for the relative $L^2$ difference between subsequent iterates to go below $10^{-3}$. The error in the reconstruction was $29\%$ and the reconstruction is plotted in Figure~\ref{fig:htrecon}. 

\begin{figure}[!ht]
\centering
\includegraphics[scale = 0.5]{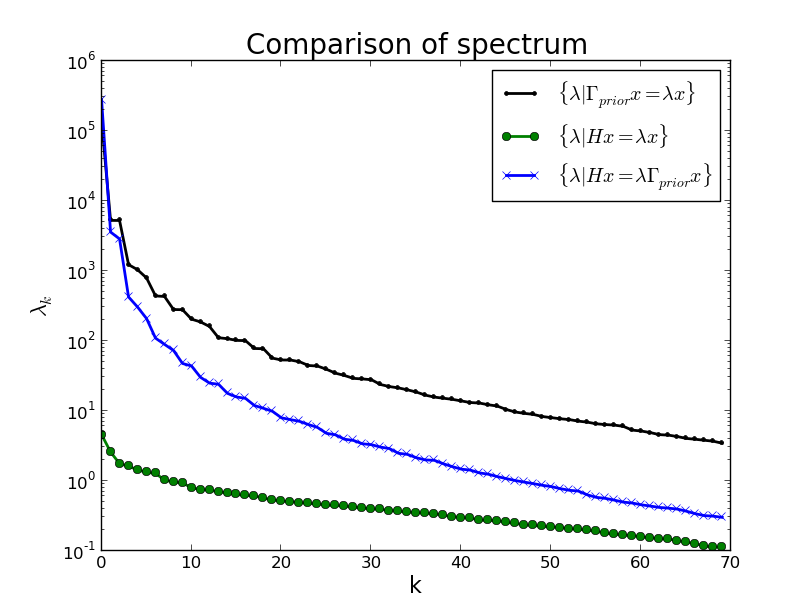}
\caption{Comparison of the eigenvalues of the data misfit Hessian $\red = J^T\noise^{-1}J$, the prior covariance matrix $\prior$ and the generalized HEP in the Figure~\ref{eqn:posteriorghep}. The Jacobian is computed by linearizing the measurement operator at the MAP estimate. As can be seen, the eigenvalues of the GHEP $\red x = \lambda \prior^{-1}x$ decays rapidly and we can retain only the largest eigenmodes above a user-defined tolerance. }
\label{fig:htcomparespectrum}
\end{figure}

\begin{figure}[!ht]
\centering
\includegraphics[scale = 0.5]{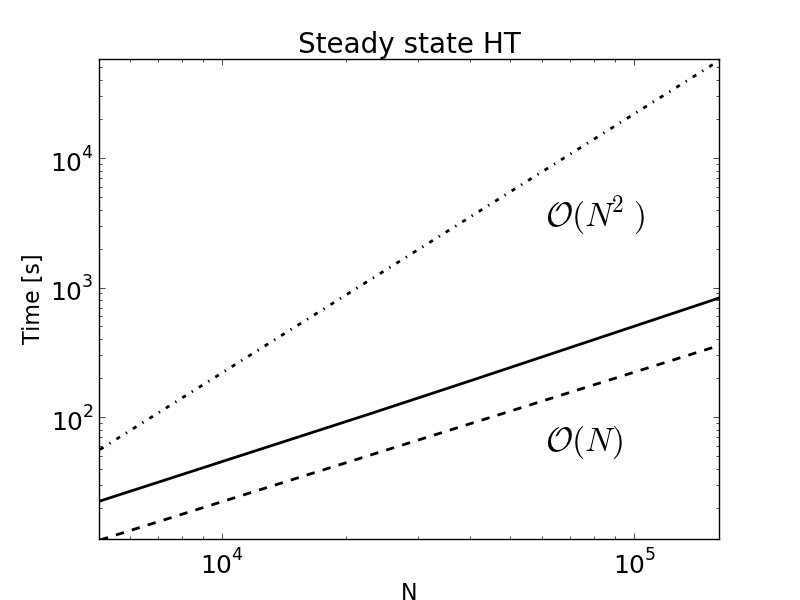}
\caption{ The CPU time (solid line) for computing the top $70$ eigenmodes of the hydraulic tomography example with unknowns discretized on grids ranging from $51\times 51$ to $401\times 401$ grid with $n_m = 1024$ measurements. The slopes of the dashed and dotted lines correspond to computational costs of $\bigO(N)$ and $\bigO(N^2)$ respectively. As can be seen, the computational cost scales linearly.It can be seen that the computational cost scales linearly with the number of unknowns.}
\label{fig:nonlinearcost}
\end{figure}

The generalized eigenvalue problem in Equation~\eqref{eqn:posteriorghep} is solved with the Jacobian computed at the MAP estimate and the eigenpairs corresponding to the $50$ largest eigenvalues are computed using the randomized GHEP algorithm described in Section~\ref{sec:randomized}. The posterior covariance matrix is computed using Equation~\eqref{eqn:postcovquasi}. The variance computed as the diagonals of the posterior covariance is  visualized in Figure~\ref{fig:htpost}. The generalized eigenvalues are plotted in Figure~\ref{fig:htcomparespectrum}. Since the true variance is hard to compute, as it requires sampling from the true posterior distribution, comparison of the accuracy against the true variance is hard to establish. However, from the analysis in Section~\ref{sec:accuracy} the computation of $F_{ss}^{-1}$ using the formula in~\eqref{eqn:fssinvapprox} a cutoff of $\lambda > 0.1 $ is a good approximation to the true $F_{ss}^{-1}$ computed at the MAP estimate. As can be seen in the Figure~\ref{fig:htcomparespectrum}, the eigenvalues of the GHEP defined in Equation~\eqref{eqn:posteriorghep} and the number of eigenvalues retained $~70$ that satisfies the cutoff $\lambda > 0.1$, is far smaller than the number of unknowns $160,801$ or the number of measurements $1024$. Finally, in Figure~\ref{fig:nonlinearcost}, we report the cost of computing the dominant eigenmodes of the GHEP defined in Equation~\eqref{eqn:posteriorghep}.  The major component in this calculation is the matvecs involving $H_\text{red}$, $\prior $ and $\prior^{-1}$. In Section~\ref{sec:posteriorapprox}, we have already argued that the costs involving the aforementioned matvecs scale linearly or almost linearly with the number of unknowns. Therefore, we expect the cost of the eigendecomposition, and as a result the variance calculations, to scale similarly with the number of unknowns. This is illustrated in Figure~\ref{fig:nonlinearcost}. As can be seen from the figure, the computational costs scale linearly with the number of unknowns.   

\section{Measures of uncertainty}\label{sec:uncertmeas}
An experimentalist interested in designing experiments to estimate parameters is interested in obtaining as much information as possible from laborious or expensive experiments. In the context of Bayesian experimental design, one tries to optimize certain criteria based on scalar measures of conditional uncertainty. Several studies have focused their attention on optimal experimental design targeted to regression-like problems. However, the large computational costs associated with storing and computing the covariance matrices prohibit their utility in optimal design for geostatistical approach. Moreover, in contrast to regression based methods the inverse problems arising from environmental sciences are highly under-determined, i.e. $n \ll m$. As a result, several uncertainty measures defined in the geostatistical context have different theoretical properties and significance compared to the (typically) over-determined regression-based problems. For a good review of uncertainty measures in the context of the geostatistical approach, please refer to~\cite{nowak2010measures}. Several measures of uncertainty are proposed based on alphabetic optimality that is prevalent in literature on regression. Because of the high computational complexity associated with storing and computing with the posterior covariance matrices, several of the optimality criteria are either computationally intensive or completely infeasible (for example, D-optimality criterion). Several of these optimality measures can be easily computed when the number of measurements are small, even when the unknowns are discretized on a very fine grid. However, as the number of measurements grow, the computational cost associated with estimating the optimality criteria can grow as either $\bigO(nm +n^3)$ in certain cases. 

 The optimality criteria are then based upon the invariants of the posterior covariance matrix, such as the trace, log determinant etc. When the posterior PDF is Gaussian, the posterior covariance matrix is provided by Equation~\eqref{eqn:posteriorlinear}. When the posterior distribution is strongly non-Gaussian, the approximation using~\eqref{eqn:postcovquasi} is no longer a good approximation to the posterior covariance and one has to resort to Monte-Carlo methods for computing these optimality criteria~\cite{huan2012simulation}. Our approach is to approximate the posterior covariance  by computing the dominant eigenmodes of the GHEPin Equation~\eqref{eqn:posteriorghep} and using the approximation described in~\eqref{eqn:fisherinverse} and the procedure described in Section~\ref{sec:posteriorapprox}. Here, we describe the calculation of a few measures of uncertainty based on the posterior covariance matrix and efficient ways to calculate it. These notations are also described in~\cite{nowak2010measures}. 

\begin{enumerate}
\item A-optimality 
\begin{equation}\label{eqn:Acriterion}
\phi_A \define \frac{1}{m+n_p}\trace\left( A \post\right) 
\end{equation} 

For $A=I$ the A-optimality criterion takes the special form $\phi_I=\frac{1}{m+n_p}\trace (\post) $.  To compute $\phi_I$, we start with the observation that 
\[ \trace (\post) = \trace(F_{ss}^{-1}) + \trace (GS^{-1}G^T) + \trace (S^{-1})\]
 Since $S$ is of size $n_p\times n_p$ we can compute $S^{-1}$ and its trace explicitly in $\bigO(n_p^3)$.  Since $G$ is $m\times n_p$ , $\trace (GS^{-1}G^T)$ can be computed in $\bigO (mn_p)$. Finally, $\trace ( F_{ss}^{-1})$ can be computed as 
\begin{equation}
\trace (F_{ss}^{-1}) \approx \trace (\prior) - \trace (U_kD_kU_k^T)
\end{equation}
which can be computed in $\bigO (mk)$ and $\trace (\prior)$ is known from the covariance kernel. In summary, the calculation of $\trace (\post)$ is $\bigO (k+n_p)m$. In practice, $n_p$ is $\bigO(1)$, so the terms involving it have negligible contribution to the overall computational cost.

For $A\neq I $, $\phi_A$ can be approximated using the Hutchinson trace estimator. Essentially,

\[ \trace (A\post) \approx \frac{1}{s} \sum_{i=1}^s v_i^T A\post v_i\] 
where the entries of $v_i$ are chosen as $\pm 1$ with equal probability. This forms an unbiased estimator for the trace of $A\post$. However, the variance of the estimator can be large. To remedy this, other sophisticated sampling schemes for estimating the trace of a matrix, based only on the availability of fast matvecs, have been proposed in~\cite{avron2011randomized}.  
\item C-optimality

 \begin{equation} \label{eqn:ccriterion}
\phi_C \define \frac{1}{m+n_p}c^T\post c 
\end{equation} 
where $c = [\partial z /\partial s^T \partial z /\partial \beta^T]^T $ is the sensitivity of a scalar model prediction $z$ with respect to the estimated parameters, yielding the prediction variance for $z$. This criterion can be easily extended to vector valued model predictions as well. For details, see~\cite{nowak2010measures}. $\phi_C$ is a special case of $\phi_A$. In particular,
\[ \phi_C = \frac{1}{m+n_p}\trace ( cc^T\post) =  \phi_{A = cc^T}\]
Given the vector $c$, it is easy to see that $\phi_C$ can be calculated in $\bigO (m\log m + km)$. Furthermore, $c$ can computed efficiently using adjoint-based techniques~\cite{nowak2010measures}.

\item  D-optimality 
\[ \phi_D  \define \logdet\left(\post\right)\]
For numerical reasons, the logarithm is considered here instead of the determinant itself. The calculation of the determinant is computationally infeasible. However, since $\post$  is a symmetric positive definite matrix 
\[ \log \text{det} \left(\post\right) = \trace \left(\log \post \right)\] 
Then, in a procedure similar to the $A$-optimality criterion, we can compute the trace of the matrix using the Hutchinson trace estimator
\[ \trace \left(\log \post \right) \approx \frac{1}{s}\sum_{i=1}^s v_i^T\log \post v_i\]
As before,  the entries of $v_i$ are chosen as $\pm 1$ with equal probability. Calculating the entries of $\log \post$ is also computationally infeasible. However, matvecs of the form $\log (\post)v_i $ can be formed efficiently using matrix-free techniques described in~\cite{chen2011computing,hale2008computing}. When the posterior PDF is Gaussian, i.e. for linear inverse problems, the D-optimality criterion is closely related to maximizing the entropy of the random variable. The entropy of a random variable $X$ with PDF $p(X)$ is defined as $H[X] =  E[ -\log p(X)] $, where $E[\cdot]$ is the expectation. For Gaussian distributions $X \sim \normal(\mu,\Sigma)$, the entropy can be calculated as $H[X] = \frac{1}{2}\log 2\pi e + \frac{1}{2} \log \det (\Sigma)$. Therefore, $\phi_D =  2 H[ s,\beta | y] - \log 2\pi e $.

An alternative way to compute D-optimality is to use the properties of determinants - the determinant of the product of matrices is the product of the determinants. Applying this result to Equation~\eqref{eqn:fisherinverse}, we have that 
\begin{align}
 \logdet(\post) = & \quad  \logdet(F_{ss}^{-1}) + \logdet (S^{-1})  \\ \nonumber
	\approx & \quad  \logdet (\prior) +  \logdet(I-D_k) + \logdet (S^{-1}) 
\end{align}
This result follows because 
\begin{align*}
 \logdet (F_{ss}^{-1}) \approx & \quad \logdet (\prior) + \logdet\left(I - \prior^{-1} U_kD_kU_k^T\right) \\
= & \quad \logdet (\prior) + \logdet \left( I - U_k^T\prior^{-1}U_kD_k \right) \\ \nonumber  
= & \quad \logdet (\prior) + \logdet (I-D_k)   
\end{align*}
where the equality in the second line followed from a straightforward application of the Sylvester's determinant lemma followed by the use of the fact that the generalized eigenvectors $U_k$ were $\prior^{-1}$-orthonormal. Further, we can treat $\logdet (\prior)$ as constant because it does not change with the measurements. 
\item E-optimality 
\[ \phi_E \define \max \text{eig} (\post)\] 
This is related to the C-optimality as $\max_{\normtwo{c} = 1} \phi_C$. The resulting vector $c$ that maximizes the C-optimality criterion is the direction of the largest variation. Computing the largest eigenvalue (and eigenvector) can be easily computed by using a few matvecs involving $\post$ using few iterations of either the power iteration or Lanczos method. It can also be efficiently estimated using the randomized GHEP algorithm that we describe in Section~\ref{sec:randomized} with $A=\post$ and $B=I$.

\end{enumerate}

\begin{figure}[!ht]
\centering
\includegraphics[scale = 0.5]{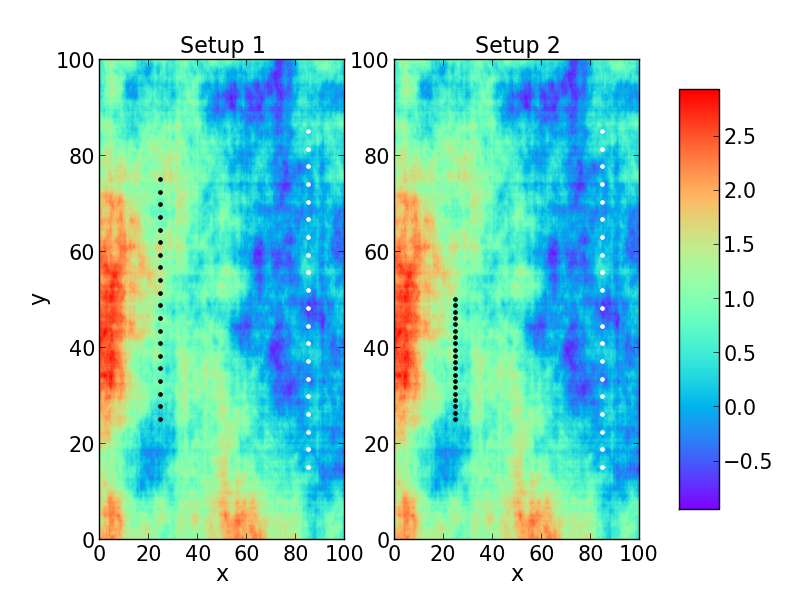}
\caption{Two different setups for comparing uncertainty measures defined in~\ref{sec:uncertmeas}. Black dots correspond to sources and white dots correspond to receivers, of which there are 20 each. The true field is discretized on a $512\times 512$ grid and is a realization of the covariance kernel given in~\eqref{eqn:raycovariance}. }
\label{fig:uncertmeas}
\end{figure}

To illustrate these ideas, we consider again the ray tomography example in Section~\ref{sec:raytomog}.
We choose two different experimental setups to compare the uncertainty measures described in Section~\ref{sec:uncertmeas}. Both setups have the same number of sources $20$ and same number of receivers $20$ discretized on a domain of size $512 \times 512$. However, they differ in the locations of the sources and as a result produce different reconstructions of the same field, see Figure~\ref{fig:uncertmeas}. But a natural question to ask is: which experimental setup is better? Intuitively, one would expect that when the sources cover a larger area the reconstruction would be better over a larger area and hence, the resulting uncertainty (measured as a function of the posterior covariance) would be lower in the entire domain. Indeed, this is the case. Quantitatively, the variance of setup $2$ is larger than that of setup $1$ and this is reflected in all the measures of uncertainty. This is listed in Table~\ref{tab:uncertmeas}. 

\begin{table}[!ht]
\centering
\begin{tabular}{|c|c|c|c|c|c|}
\hline
Setup & $\phi_A$ & $\phi_C$ & $\tilde{\phi}_D$ & $\phi_E$ & $\trace{S^{-1}}$  \\ \hline
$1$ & $64.64$ & $5.22$ & $-823.18$ &$ 26.28$ & $6.31\times 10^{-4}$\\ 
$2$ & $72.33$ & $6.27$ & $-688.81$ & $26.50$ & $6.46\times 10^{-4}$\\ \hline
\end{tabular}
\caption{Comparison of the various measures of uncertainty defined in Section~\ref{sec:uncertmeas} applied to the setups described in Figure~\ref{fig:uncertmeas}. }
\label{tab:uncertmeas}
\end{table}

For the computation of $\phi_A$, we chose $A = I$, whereas for computing $\phi_C$ we chose $c$ as a vector of ones. We observe that $\tilde{\phi}_D$ is the most sensitive measure of uncertainty because the difference in $\tilde{\phi}_D$ between the two setups is the largest, even in relative magnitude. The least sensitive is the E-optimality criterion which only takes into account the largest eigenmode. $\phi_A$ measures the trace of the posterior covariance matrix. Both $\phi_A$ and $\phi_D$ account for both large-scale and small-scale variations since it takes into account all the eigenmodes of the posterior GHEP in Equation~\eqref{eqn:posteriorghep}. While computing $\tilde{\phi}_D$ it is required to compute $\logdet (I-D_k)$, and $\phi_A$ requires weighted sum of the $D_k)$. Therefore, small changes made to the lowest eigenmodes (corresponding to fine-scale variations)  have a higher impact in the calculations of $\phi_D$, rather than $\phi_A$. By the special choice of the vector $c$, which is chosen to be a vector of ones,  $\phi_C$ is proportional to the variance of predicting the global mean and can also be interpreted as proportional to the average conditional integral scale~\cite{nowak2010measures}. Since in setup 1 the sources are spread over a larger area and captures the large scale variability well but has difficulty with small scale features, leading to smaller integral scale and lower $\phi_C$. By a similar reasoning $\trace ({S}^{-1})$, which represents the uncertainty in the parameters $\beta$, is larger for setup $2$.

\section{Conclusions and Future Work}
Several papers have focused their attention on solving inverse problems. The problem of quantifying the uncertainty associated with the estimate is challenging because of the high dimensionality of the space of unknowns and the high computational costs involving the representing unknown random fields and computing the solution to the forward . We consider the problem of quantifying the associated uncertainty by computing an approximation to the posterior covariance matrix which is roughly composed of two terms - the prior covariance matrix and a low rank term that contains the dominant eigenmodes of a generalized Hermitian eigenvalue problem (GHEP) involving the  misfit portion of the Hessian and the inverse prior covariance matrix. For several inverse problems, the spectrum of the eigenproblem, described above, decays rapidly and we are justified in retaining only the largest eigenvalues satisfying a cutoff $\lambda > \varepsilon$. We have provided an efficient algorithm for computing this low-rank representation using a randomized algorithm that deliberately avoids expensive computations involving square-root (or its inverse) of the prior covariance matrix. The accuracy of the low-rank representation and the parameters controlling the number of retained eigenvalues have been explored. Finally, we have shown how to approximately compute measures of uncertainty based on the approximate representation of the posterior covariance matrix. The resulting algorithms are highly scalable and their performance has been demonstrated on two challenging applications - ray based tomography and steady-state hydraulic tomography. 

For nonlinear problems we have used the approximate posterior distribution which is approximated to be a Gaussian by linearizing the measurement operator at the MAP point. However, for highly nonlinear problems this approximation may not be reasonable. In such cases, one may have to resort to sampling from the full posterior in order to quantify the uncertainty~\cite{martin2012stochastic}. However, such a process would be challenging for much the same reasons - high dimensionality of the input spaces and the high computational costs associated with the forward problem. A limitation of our analysis is that it is restricted to the case for which the priors can be represented as a Gaussian random field. Future work is necessary to extend our analysis to different priors.

In future, we would like to explore the possibility of using the uncertainty measures described in Section~\ref{sec:uncertmeas} to optimize the control variables in an experimental setup in order to obtain the most amount of information from laborious or expensive experiments. Finally, another challenging application that would benefit from an efficient representation of posterior covariance, is data assimilation using Kalman filter. Implementing Kalman filter can be computationally expensive because of the high computational costs in updating the posterior distribution. For a certain class of problems, we anticipate that the posterior covariance matrix can be represented as the combination of the prior and a low-rank term, that can be recursively updated efficiently. This could lead to a highly scalable implementation of the Kalman filter. This will be discussed in subsequent papers. 

\section{Acknowledgements}

The major work on this paper was completed when the first author was a PhD candidate in ICME Stanford University and he would like to thank Ivy Huang for all her help and support during this process and beyond.  The authors were supported by NSF Award 0934596, Subsurface Imaging and Uncertainty Quantification.  We thank Tania Bakhos for a careful reading of the manuscript and providing valuable suggestions to improve readability.

\bibliographystyle{plain}
\bibliography{main}

\end{document}